\documentclass[manmat,numbook,envcountsame]{svjour}

\usepackage{amsmath,latexsym,amssymb,times,mathptm}

\def\sqr#1#2{{\vcenter{\hrule height.#2pt
        \hbox{\vrule width.#2pt height#1pt \kern#1pt
                \vrule width.#2pt}
        \hrule height.#2pt}}}
\def\square{\mathchoice\sqr64\sqr64\sqr{4}3\sqr{3}3}
\def\QED{\hfill$\square$}

\def\tratto{\mbox{\rule{2mm}{.2mm}$\;\!$}}

\begin{document}

\title{Core of projective dimension one modules\thanks{The first author
was partially supported by a Faculty Research Fellowship from the
University of Kentucky. The second and third authors gratefully
acknowledge partial support from the NSF.}}

\titlerunning{Core of projective dimension one modules}
\authorrunning{A. Corso et al.}

\author{Alberto Corso \and Claudia Polini \and Bernd Ulrich}

\institute{{\sc A. Corso} \\
Department of Mathematics, University of Kentucky, Lexington, KY
40506, USA \\ (e-mail$\colon$ corso@ms.uky.edu) \and
{\sc C. Polini} \\
Department of Mathematics, University of Notre Dame, Notre Dame,
IN 46556, USA \\(e-mail$\colon$ cpolini@nd.edu) \and
{\sc B. Ulrich} \\
Department of Mathematics, Purdue University, West Lafayette, IN 47907, USA \\
(e-mail$\colon$ ulrich@math.purdue.edu)}

\date{Received: 10 June 2002/}

\maketitle

\begin{abstract}
The core of a projective dimension one module is computed
explicitly in terms of Fitting ideals. In particular, our formula
recovers previous work by R. Mohan on integrally closed
torsionfree modules over a two-dimensional regular local ring.

\medskip

\noindent \subclassname{Primary 13C10; Secondary 13A30, 13B21.}
\end{abstract}

\section{Introduction}

There is an extensive literature on the {\it Rees algebra}
${\mathcal R}(I)$ of an $R$-ideal $I$ in a Noetherian ring $R$.
One reason for that is the major role that ${\mathcal R}(I)$ plays
in Commutative Algebra and Algebraic Geometry, since ${\rm
Proj}({\mathcal R}(I))$ is the blowup of ${\rm Spec}(R)$ along the
subscheme defined by $I$.

There are several reasons for studying Rees algebras of finitely
generated $R$-modules $E$ as well. For instance, Rees algebras of
modules include the so called multi-Rees algebras, which
correspond to the case of direct sums of ideals. Furthermore,
symmetric algebras of modules arise naturally as coordinate rings
of certain correspondences in Algebraic Geometry, and the
projection of these varieties often requires the killing of
torsion. This takes us back to the study of  Rees algebras of
modules: We refer to the articles by D. Eisenbud, C. Huneke and B.
Ulrich \cite{EHU} and by A. Simis, B. Ulrich and W.V. Vasconcelos
\cite{RAM} for additional motivation as well as a rich list of
references. We stress, though, that this is not a routine
generalization of what happens for ideals. This is essentially due
to the fact that the remarkable interaction that exists between
the Rees algebra and the associated graded ring of an ideal is
missing in the module case.

The importance of {\it reductions} in the study of Rees algebras
of ideals has long been noticed by commutative algebraists.
Roughly, a reduction is a simplification of the given ideal which
carries the most relevant information about the original ideal
itself. In this sense, the study of the {\it core} of an ideal --
by which we mean the intersection of all the reductions of the
ideal -- helps finding uniform properties shared by all
reductions. We refer to \cite{RS,HS,CPU,CPU2} for more details as
well as the few known explicit formulas on the core of ideals.
Another motivation for studying the core of an ideal is given by
the celebrated Brian\c{c}on-Skoda theorem.

Similarly to what happens for ideals, it is useful to study Rees
algebras of modules via the ones of reduction modules. Like in the
ideal case, one then defines ${\rm core}(E)$ to be the
intersection of all reductions of $E$. In the setting of a local
Gorenstein ring $R$, our goal is to give an explicit formula for
the core in terms of a Fitting ideal of $E$, when $E$ is a
finitely generated $R$-module with rank $e$, projective dimension
one, analytic spread $\ell = \ell(E) \geq e+1$, and property
$G_{\ell-e+1}$. Interestingly enough, this is the same class of
$R$-modules that has been studied in \cite{CP3} in the context of
cancellation theorems. One of the equivalences of
Theorem~\ref{pd1} below says that
\[
{\rm core}(E) = {\rm Fitt}_{\ell}(E) \cdot E
\]
if and only if the reduction number of $E$ is at most $\ell-e$.
This result in particular covers earlier work by Mohan \cite{M}.

To prove Theorem~\ref{pd1} we reduce the rank of $E$ by factoring
out general elements. The task then becomes to assure that this
inductive procedure preserves our assumptions and conclusions.
This requires a sequence  of technical lemmas. The corresponding
results for `generic' instead of `general' elements have been
shown in \cite{RAM} -- they are useless however in our context
since they require purely transcendental residue extensions which
a priori may change the core.

\section{The main result}

We begin by reviewing some definitions and basic facts about Rees
algebras of modules. Let $R$ be a Noetherian ring and let $E$ be a
finitely generated $R$-module with rank $e>0$. The {\it Rees
algebra} ${\mathcal R}(E)$ of $E$ is the symmetric algebra of $E$
modulo its $R$-torsion. If this torsion vanishes $E$ is said to be
of {\it linear type}. Let $U \subset E$ be a submodule. One says
that $U$ is a {\it reduction} of $E$ or, equivalently, $E$ is {\it
integral} over $U$ if ${\mathcal R}(E)$ is integral over the
$R$-subalgebra generated by $U$. Alternatively, the integrality
condition is expressed by the equations ${\mathcal R}(E)_{r+1} = U
{\mathcal R}(E)_r$ with $r \gg 0$. The least integer $r \geq 0$
for which this equality holds is called the {\it reduction number
of $E$ with respect to $U$} and denoted by $r_U(E)$. Observe that
if $E$ is of linear type then $E$ has no proper reductions. For
any reduction $U$ of $E$ the module $E/U$ is torsion, hence $U$
has a rank and ${\rm rank}\, U = {\rm rank}\, E$. If $R$ is
moreover local with residue field $k$ then the Krull dimension of
${\mathcal R}(E)\otimes_R k$ is called the {\it analytic spread}
of $E$ and is denoted by $\ell(E)$. A reduction of $E$ is said to
be {\it minimal} in case it is minimal with respect to inclusion.
If $k$ is infinite then minimal reductions of $E$ always exist and
their minimal number of generators equals $\ell(E)$. The {\it
reduction number} $r(E)$ of $E$ is defined to be the minimum of
$r_U(E)$, where $U$ ranges over all minimal reductions of $E$. The
{\it core of $E$}, denoted ${\rm core}(E)$, is the intersection of
all reductions $($minimal reductions if $R$ is local and $k$ is
infinite$)$ of $E$. We write ${\rm Fitt}_i(E)$ for the $i^{\rm
th}$ Fitting ideal of $E$. Finally, we recall that $E$ is said to
satisfy condition $G_s$, for an integer $s\geq 1$, if $\mu
(E_{\mathfrak p})\leq \dim R_{\mathfrak p} +e-1$ whenever $1\leq
\dim R_{\mathfrak p}\leq s-1$.

\medskip

Our first lemma establishes the existence of `superficial' elements for
modules:

\begin{lemma}\label{5.1}
Let $R$ be a local Cohen-Macaulay ring with infinite residue
field, let $E$ be a finitely generated $R$-module having rank
$\geq 2$, and let $U$ be a reduction of $E$. If $x$ is a general
element of $U$ then $Rx \simeq R$ and $x$ is regular on ${\mathcal
R}(E)$. Writing $\overline{E} = E/Rx$ and $K$ for the kernel of
the natural epimorphism from $\overline{\mathcal R} = {\mathcal
R}(E)/(x)$ onto ${\mathcal R}(\overline{E})$, we have that $K =
H^0_{\overline{\mathcal R}_{+}}(\overline{\mathcal R})$, in
particular $K_n=0$ for $n \gg 0$.
\end{lemma}
\begin{proof}
Since ${\rm rank}\, U = {\rm rank}\, E \geq 2$ and $x \in U$ is
general, we have that $x$ is basic for $U$ locally in codimension
one \cite[A]{EE}. Hence $R x \simeq R$ and $x$ generates a free
summand of $E$ locally at every minimal prime of $R$. It follows
that $x$ is regular on ${\mathcal R}(E)$, $K$ is the $R$-torsion
of $\overline{\mathcal R}$, and $H^0_{\overline{\mathcal
R}_{+}}(\overline{\mathcal R}) \subset K$.

To prove the inclusion $K \subset H^0_{\overline{\mathcal R}_{+}}
(\overline{\mathcal R})$, write ${\mathcal R} = {\mathcal R}(E)$,
$A = \overline{\mathcal R}$, $B = A / H_{A_{+}}^0(A)$. We need to
show that $B$ is $R$-torsionfree. We first prove that the ring $B$
satisfies $S_1$, or equivalently that every associated prime of
the ring $A$ not containing $A_{+}$ is minimal. To this end we
consider the following subsets of ${\rm Spec}({\mathcal R})$,
${\mathcal P} = \{ P \, | \,  {\rm dim}\, {\mathcal R}_{\ P} \geq
2 > {\rm depth}\, {\mathcal R}_{\ P} \} \setminus V({\mathcal
R}_{\ +})$ and $V(L)= \{ P \, | \, {\mathcal R}_{\ P} \ {\rm is \
not \ } S_2 \}$ for $L$ some ${\mathcal R}$-ideal, which exists by
\cite[6.11.8 and 6.11.9$({\it i})$]{EGA}. As ${\mathcal R}$
satisfies $S_1$, $L$ contains an ${\mathcal R}$-regular element
$y$. Now ${\mathcal P} \subset {\rm Ass}_{\mathcal R}({\mathcal
R}/(y))$, showing that ${\mathcal P}$ is a finite set $($see also
\cite[3.2]{F}$)$. As $V({\mathcal R}_{\ +}) = V(U {\mathcal R})$,
a general element $x \in U$ is not contained in any prime of
${\mathcal P}$. Therefore every associated prime of $A$ not
containing $A_{+}$ is indeed minimal.

It remains to show that every minimal prime $P$ of $A$ contracts
to a minimal prime of $R$. To this end let $Q$ be the preimage of
$P$ in ${\mathcal R}$ and  write ${\mathfrak p} = Q \cap R = P
\cap R$. Set $e = {\rm rank}\, E$ and consider the subset
${\mathcal F} = {\rm Min}({\rm Fitt}_e(E) {\mathcal R}) \setminus
V({\mathcal R}_{\ +})$ of ${\rm Spec}({\mathcal R})$. Again since
$V({\mathcal R}_{\ +}) = V(U{\mathcal R})$, a general element $x
\in U$ is not contained in any prime of ${\mathcal F}$. Thus since
${\rm ht}\, {\mathcal R}_{\ +} = e \geq 2 > 1 = {\rm ht}\, Q$ and
${\rm ht}\, {\rm Fitt}_e(E) {\mathcal R} \geq 1$, $Q$ cannot
contain ${\rm Fitt}_e(E) {\mathcal R}$. Therefore $E_{\mathfrak
p}$ is free, hence $E_{\mathfrak p} =U_{\mathfrak p}$ and ${\rm
dim}\, R_{\mathfrak p} \leq {\rm dim}\, {\mathcal R}_{\ Q} = 1$.
Thus $x$ generates a free summand of $E_{\mathfrak p}
=U_{\mathfrak p}$ , and so $A_{\mathfrak p} = {\mathcal R}_{\
\mathfrak p}/(x)$ is a polynomial ring over $R_{\mathfrak p}$.
Therefore ${\rm dim}\, R_{\mathfrak p} \leq {\rm dim}\, A_P = 0$.
\QED
\end{proof}

The next lemma is the main ingredient to set up the inductive
procedure in the proof of Theorem~\ref{pd1}. We show how the
properties of a module $E$ change after factoring out a general
element.

\begin{lemma}\label{5.2}
In addition to the assumptions and notation of\/ {\rm
Lemma~\ref{5.1}} write $e = {\rm rank}\, E$ and let {\rm
`${}^{\tratto}$'} denote images in $\overline{E}$.
\begin{description}[(a)]
\item[$(${\it a}$)$]
$\overline{E}$ has rank $e-1$;

\item[$(${\it b}$)$]
if\/ $E$ is torsionfree and satisfies $G_s$ and if\/ ${\rm ht}\, U
:_R E \geq s$ for some $s \geq 2$, then $\overline{E}$ is
torsionfree and satisfies $G_s$;

\item[$(${\it c}$)$]
if\/ $V \subset E$ is a submodule containing $x$ and
$\overline{V}$ is a reduction of $\overline{E}$, then $V$ is a
reduction of $E$;

\item[$(${\it d}$)$]
$\ell(\overline{E}) = \ell(E)-1$;

\item[$(${\it e}$)$]
if\/ $U$ is a minimal reduction of $E$ then $\overline{U}$ is a
minimal reduction of $\overline{E}$ and $r(\overline{E}) \leq
r_U(E)$;

\item[$($\/{\it f}\/$)$]
if\/ ${\mathcal R}(E)$ or ${\mathcal R}(\overline{E})$ satisfies
$S_2$ then $K=0$.
\end{description}
\end{lemma}
\begin{proof}
Let $k$ denote the residue field of $R$. Part $(${\it a}$)$
follows since $Rx \simeq R$ by Lemma~\ref{5.1}. As for $(${\it
b}$)$, notice that $U$ satisfies $G_s$ and then $\overline{U}$ has
the same property by \cite[3.5 and its proof]{MY}. Thus
$\overline{E}$ is $G_s$ since ${\rm ht}\, \overline{U} \colon
\overline{E} \geq s$. It now follows that $\overline{E}$ is
torsionfree, because $E$ is torsionfree and $\overline{E}$ is free
in codimension $1$. We now prove $(${\it c}$)$. By
Lemma~\ref{5.1}, the kernel of the natural map from $A =
({\mathcal R}(E)/V{\mathcal R}(E)) \otimes_R k$ to $B = ({\mathcal
R}(\overline{E})/\overline{V}{\mathcal R}(\overline{E})) \otimes_R
k$ is contained in $H_{A_{+}}^0(A)$. As ${\rm dim}\, B = 0$ it
follows that ${\rm dim}\, A = 0$, hence $V$ is a reduction of $E$.
To prove $(${\it d}$)$, we may assume that the image of $x$ is
part of a system of parameters of ${\mathcal R}(E) \otimes_R k$.
Thus $\ell(\overline{E}) \leq \ell(E)-1$. Now part $(${\it c}$)$
gives the asserted equality. As for $(${\it e}$)$, we may choose
$x$ to be part of a minimal generating set of $U$. Thus $(${\it
d}$)$ implies that $\overline{U}$ is a minimal reduction of
$\overline{E}$, and hence $r(\overline{E}) \leq
r_{\overline{U}}(\overline{E}) \leq r_U(E)$. Finally $($\/{\it
f}\/$)$ follows from Lemma~\ref{5.1}. This is obvious in case
${\mathcal R}(E)$ is $S_2$. If on the other hand ${\mathcal
R}(\overline{E})$ satisfies $S_2$, one argues as in \cite[the
proof of 3.7]{RAM}. \QED
\end{proof}

\begin{lemma}\label{5.3}
Let $R$ be a local Cohen-Macaulay ring with infinite residue
field, and let $E$ be a finitely generated $R$-module with ${\rm
proj\, dim\,}E = 1$. Write $e = {\rm rank\,} E$, $\ell = \ell(E)$,
and assume that $E$ satisfies $G_{\ell-e+1}$ and is torsionfree
locally in codimension $1$. If\/ $U$ is any minimal reduction of
$E$ then $U/(U :_R E)U \simeq (R/U :_R E)^{\ell}$.
\end{lemma}
\begin{proof}
We use induction on $e \geq 1$. Notice that $E$ is torsionfree. If
$e=1$ then $E$ is isomorphic to a perfect ideal of grade $2$. Now
the assertion follows from \cite[2.5]{CPU2}. Thus we may assume $e
\geq 2$. Let $x$ be a general element of $U$ and let
`${}^{\tratto}$' denote images in $\overline{E} = E/Rx$. By
Lemma~\ref{5.1}, $Rx \simeq R$. Notice that $E$ is free locally in
codimension $1$, thus ${\rm ht}\, U \colon E \geq 2$. Furthermore
according to \cite[Proposition 4]{Avr}, $E$ is of linear type
locally in codimension $\ell-e$, hence ${\rm ht}\, U \colon E \geq
\ell-e+1$. Therefore by Lemma~\ref{5.2}$(${\it a}$)$,$(${\it
b}$)$,$(${\it d}$)$,$(${\it e}$)$, $\overline{E}$ satisfies the
same assumptions as $E$ with ${\rm rank}\, \overline{E} = e-1$ and
$\ell(\overline{E}) = \ell-1$, and $\overline{U}$ is a minimal
reduction of $\overline{E}$. These facts remain true if we replace
$x$ by $x_i$, for a suitable generating set $\{ x_1, \ldots,
x_{\ell} \}$ of $U$. Now the asserted isomorphism holds if and
only if ${\rm Fitt}_{\ell-1}(U) \subset U \colon E$, or
equivalently $(Rx_1 + \ldots + \widehat{Rx_i}+\ldots+Rx_{\ell})
\colon x_i \subset U \colon E$ for every $1 \leq i \leq \ell$.
Taking $x \in \{ x_1, \ldots, \widehat{x}_i, \ldots, x_{\ell} \}$,
we have $(Rx_1 + \ldots + \widehat{Rx_i}+\ldots+Rx_{\ell}) \colon
x_i = \overline{(Rx_1 + \ldots + \widehat{Rx_i}+\ldots+Rx_{\ell})}
\colon \overline{x}_i \subset \overline{U} \colon \overline{E} = U
\colon E$, where the middle inclusion holds by induction
hypothesis. \QED
\end{proof}

We are now ready to prove our main result, which characterizes the
shape of the core of a module with projective dimension one in
terms of conditions either on the reduction number or on a Fitting
ideal of the module. The corresponding result for ideals has been
shown in \cite[2.6 and 3.4]{CPU2}.

\begin{theorem}\label{pd1}
Let $R$ be a local Gorenstein ring with infinite residue field and let
$E$ be a finitely generated $R$-module with ${\rm proj\, dim\,}E = 1$.
Write $e = {\rm rank\,} E$, $\ell=\ell(E)$, and assume that $E$ satisfies
$G_{\ell-e+1}$ and is torsionfree locally in codimension $1$.
The following conditions are equivalent:
\begin{description}[(a)]
\item[$(${\it a}$)$]
$(U :_R E) E \subset {\rm core}(E)$ for some minimal reduction $U$
of $E$;

\item[$(${\it b}$)$]
$(U :_R E) U = (U :_R E) E = {\rm core}(E)$ for every minimal
reduction $U$ of $E$;

\item[$(${\it c}$)$]
${\rm core}(E) = {\rm Fitt}_{\ell}(E) \cdot E$;

\item[$(${\it d}$)$]
$U :_R E$ does not depend on the minimal reduction $U$ of $E$;

\item[$(${\it e}$)$]
$U :_R E = {\rm Fitt}_{\ell}(E)$ for every minimal reduction $U$ of $E$;

\item[$($\/{\it f}\/$)$]
the reduction number of $E$ is at most $\ell-e$.
\end{description}
\end{theorem}
\begin{proof}
Notice that $E$ is free locally in codimension $1$ and
torsionfree. Furthermore $\ell \geq e+1$ by \cite[4.1$(${\it
a}$)$]{RAM}.  Let $U$ be any minimal reduction of $E$. By
\cite[Proposition 4]{Avr}, $E$ is of linear type locally in
codimension $\ell-e$, hence ${\rm ht}\, U \colon E \geq \ell-e+1$.
Thus $U \colon E = {\rm Fitt}_0(E/U)$ according to
\cite[3.1$($2$)$]{BE}. In particular ${\rm Fitt}_{\ell}(E) = \sum
(U \colon E)$ with $U$ ranging over all minimal reductions of $E$.
This establishes the equivalence $(${\it d}$)$ $\Leftrightarrow$
$(${\it e}$)$ and the implications $(${\it b}$)$ $\Rightarrow$
$(${\it c}$)$ $\Rightarrow$ $(${\it a}$)$. The implication $(${\it
d}$)$ $\Rightarrow$ $(${\it a}$)$ on the other hand is obvious.

Now we are going to prove the remaining equivalences by induction
on $e \geq 1$. If $e=1$ then $E$ is isomorphic to a perfect ideal
of grade $2$, and the theorem follows from \cite[2.6$($3$)$]{CPU2}
and \cite[3.4]{CPU2}. Thus we may assume that $e \geq 2$. Let $x$
be a general element of $U$ and let `${}^{\tratto}$' denote images
in $\overline{E} = E/Rx$. By Lemma~\ref{5.1}, $Rx \simeq R$ and
$x$ is regular on ${\mathcal R}(E)$. Furthermore
Lemma~\ref{5.2}$(${\it a}$)$,$(${\it b}$)$,$(${\it d}$)$,$(${\it
e}$)$ shows that $\overline{E}$ satisfies the same assumptions as
$E$ with ${\rm rank}\, \overline{E} = e-1$, $\ell(\overline{E}) =
\ell-1$ and $\overline{U}$ a minimal reduction of $\overline{E}$.
Finally $\overline{{\rm core}(E)} \subset {\rm
core}(\overline{E})$ by Lemma~\ref{5.2}$({\it c})$. These facts
remain true if we replace $x$ by $x_i$, for a suitable generating
set $\{ x_1, \ldots, x_{\ell} \}$ of $U$.

$(${\it a}$)$ $\Rightarrow$ $(${\it f}\/$)$: Take $U$ as in
$(${\it a}$)$. Since $\overline{{\rm core}(E)} \subset {\rm
core}(\overline{E})$, part $(${\it a}$)$ gives $(\overline{U}
\colon \overline{E})\overline{E} \subset {\rm
core}(\overline{E})$. Now the induction hypothesis implies that
$r(\overline{E}) \leq \ell(\overline{E})-(e-1)$. Hence ${\mathcal
R}(\overline{E})$ is Cohen-Macaulay by \cite[4.7$(${\it
a}$)$]{RAM}. Thus Lemma~\ref{5.2}$($\/{\it f}\/$)$ shows that
${\mathcal R}(E)$ is Cohen-Macaulay, and then $r(E) \leq \ell-e$
again by \cite[4.7$(${\it a}$)$]{RAM}.

$($\/{\it f}\/$)$ $\Rightarrow$ $(${\it e}$)$ and $(${\it b}$)$:
If $($\/{\it f}\/$)$ holds then ${\mathcal R}(E)$ is
Cohen-Macaulay by \cite[4.7$(${\it a}$)$]{RAM}.
Lemma~\ref{5.2}$($\/{\it f}\/$)$ then shows that ${\mathcal
R}(\overline{E}) \simeq {\mathcal R}(E)/(x)$ is Cohen-Macaulay as
well. Now again according to \cite[4.7$(${\it a}$)$]{RAM},
$r(\overline{E}) \leq \ell(\overline{E})-(e-1)$. To establish
$(${\it e}$)$ we use the isomorphism $\omega_{{\mathcal R}(E)}
\simeq {\rm Fitt}_{\ell}(E) {\mathcal R}(E)[-e]$ proved in
\cite[4.10]{RAM}. From the corresponding formula for
$\omega_{{\mathcal R}(\overline{E})}$  and the isomorphism
$\omega_{{\mathcal R}(E)} \otimes_{{\mathcal R}(E)} {\mathcal
R}(E)/(x) \simeq \omega_{{\mathcal R}(\overline{E})}[-1]$, we
deduce that ${\rm Fitt}_{\ell}(E) = {\rm
Fitt}_{\ell-1}(\overline{E})$. On the other hand by induction
hypothesis, ${\rm Fitt}_{\ell-1}(\overline{E}) = \overline{U}
\colon \overline{E}$. Since $\overline{U} \colon \overline{E} = U
\colon E$ we obtain the equality ${\rm Fitt}_{\ell}(E) = U \colon
E$ asserted in $(${\it e}$)$. Now $(${\it e}$)$ being established,
part $(${\it d}$)$  holds as well and hence $(U \colon E)E \subset
{\rm core}(E)$ for every minimal reduction $U$ of $E$. Thus to
prove $(${\it b}$)$ it suffices to show that ${\rm core}(E)
\subset (U \colon E)U$. By the induction hypothesis, ${\rm
core}(\overline{E}) = (\overline{U} \colon
\overline{E})\overline{U}$. Since $\overline{{\rm core}(E)}
\subset {\rm core}(\overline{E})$ we deduce that ${\rm core}(E)
\subset Rx + (U \colon E)U$. Replacing $x$ by $x_i$ we obtain
\[
{\rm core}(E) \subset \bigcap_{i=1}^{\ell} (Rx_i+(U \colon E)U).
\]
However the latter module is $(U \colon E)U$. This is obvious if
$U \colon E = R$ and follows from Lemma~\ref{5.3} if $U \colon E
\not= R$ since then the images of $x_1, \ldots, x_{\ell}$ form a
basis of the free module $U/(U \colon E)U$ over the ring $R/U
\colon E$. \QED
\end{proof}

Our first corollary says that in the setting of Theorem~\ref{pd1},
formation of the core is indeed compatible with flat local base
change.

\begin{corollary}
Let $R$ be a Noetherian local ring with infinite residue field and
let $E$ be a finitely generated $R$-module with ${\rm proj\,
dim}\, E =1$. Write $e={\rm rank}\, E$, $\ell=\ell(E)$, and assume
that $E$ satisfies $G_{\ell-e+1}$ and is torsionfree locally in
codimension $1$ and that $r(E) \leq \ell-e$. For any flat local
homomorphism $R \longrightarrow S$ with $S$ a local Gorenstein
ring, one has ${\rm core}(E \otimes_R S) =({\rm core}(E)) S$.
\end{corollary}

Our next corollary deals with a formula for the core of modules
with projective dimension one that are presented by a matrix with
linear entries.

\begin{corollary}
Let $R = k[x_1, \ldots, x_d]_{(x_1, \ldots, x_d)}$ be the localization of
a polynomial ring in $d \geq 2$ variables over an
infinite field, let ${\mathfrak m}$ be the maximal
ideal of $R$, and let $E$ be a finitely generated $R$-module with
${\rm proj\, dim}\, E =1$, presented by a matrix whose entries are
linear forms. Write $e = {\rm rank}\, E$, $n = \mu(E)$, and assume
that $E$ satisfies $G_d$. One has
\[
{\rm core}(E) = {\mathfrak m}^{n-e-d+1} E.
\]
\end{corollary}
\begin{proof}
We may assume that $n > d+e-1$, since otherwise $E$ is of linear
type according to \cite[4.11]{RAM}. But then by the same result,
$\ell=\ell(E)=d+e-1$. Furthermore ${\mathcal R}(E)$ is
Cohen-Macaulay. Thus \cite[4.7$(${\it a}$)$]{RAM} shows that $r(E)
\leq \ell-e$ and ${\rm Fitt}_{\ell}(E) = {\rm Fitt}_0(E/U)$ for
some minimal reduction $U$ of $E$. Now according to
Theorem~\ref{pd1}, ${\rm core}(E)={\rm Fitt}_0(E/U) \cdot E$. On
the other hand, ${\rm Fitt}_0(E/U)$ has height $d$ and is the
ideal of maximal minors of an $n-e-d+1$ by $n-e$ matrix with
linear entries. Thus ${\rm Fitt}_0(E/U) = {\mathfrak
m}^{n-e-d+1}$. \QED
\end{proof}

Theorem~\ref{pd1} also recovers the next result of Mohan about the
core of integrally closed modules over a two-dimensional regular
local ring \cite[2.10]{M}.

\begin{corollary}
Let $R$ be a two-dimensional regular local ring with infinite
residue field and let $E$ be a finitely generated torsionfree
$R$-module of rank $e$. If $E$ is integrally closed then
\[
{\rm core}(E) = {\rm Fitt}_{e+1}(E) \cdot E.
\]
\end{corollary}
\begin{proof}
We may assume that $E$ is not free, in which case $\ell(E)=e+1$ by
\cite[4.1$(${\it a}$)$]{RAM}. Since $r(E) \leq 1$ according to
\cite[4.1$(${\it i}$)$]{KK}, the assertion now follows from
Theorem~\ref{pd1}. \QED
\end{proof}

\end{document}